\documentclass{amsart}
\usepackage{cancel} 
\usepackage{layout} 
\usepackage{amsmath,amssymb,amsfonts,mathrsfs,amscd,euscript,enumerate,calc,verbatim}
\usepackage{float}

\usepackage{latexsym}

\newtheorem{theorem}{Theorem}[section]
\newtheorem{lemma}[theorem]{Lemma}

\newtheorem{conjecture}[theorem]{Conjecture}
\newtheorem{corollary}[theorem]{Corollary}

\newtheorem{definition}[theorem]{Definition}

\def\F{{\mathbb F}}
\def\U{{\mathbb U}}

\def\Fq{{\mathbb F}_q}
\def\Fqm{{\mathbb F}_{q^m}}
\def\Fqmn{{\mathbb F}_{q^{mn}}}

\def\GLm{{\mathbb {GL}}_{m}}
\def\GLmn{{\mathbb{GL}}_{mn}}

\newcommand{\TSRP}{\operatorname{TSRP}}
\newcommand{\TSR}{\operatorname{TSR}}

\usepackage{url}

\begin{document}

\title[Primitive Transformation shift registers over finite fields]
{Primitive Transformation shift registers over finite fields.}
\author{Ambrish Awasthi}
\address{Department of Mathematics, \newline \indent IIT Delhi, Hauz Khas, New Delhi 110016, India}
\email{ambrishawasthi@yahoo.com}
\author{Rajendra K. Sharma}
\address{Department of Mathematics
\newline \indent
IIT Delhi, Hauz Khas, New Delhi 110016, India}
\email{rksharma@maths.iitd.ac.in}



\begin{abstract}
We consider the problem of existence and enumeration of primitive TSRs of order $n$ over any finite field. Here we prove the existence of primitive TSRs of order two over finite fields of characteristic $2$ and establish an equivalence between primitive TSRs and primitive polynomials of special form. A conjecture regarding the existence of these special type of primitive polynomials is submitted by us along with some experimental verification. Further we have attempted to enumerate primitive TSRs of order $2$ over finite fields of characteristic $2$. Finally we give a general search algorithm for primitive TSRs of odd order over any finite field and in particular of order two over fields of characteristic $2$.
\end{abstract}
\date{\today}
\maketitle

\section{Introduction}
\noindent
Linear feedback shift registers (LFSRs) are systems consisting of 
a homogeneous linear recurrence relation over $\Fq$. They have wide applications in cryptography
and are particularly useful for generating pseudorandom sequences in stream ciphers
refer \cite{GG,LN}. Sequences with maximum period are a necessary prerequisite for cryptographic applications.
LFSRs which generate such sequences are known as primitive LFSRs. The characteristic polynomial of such LFSRs are primitive in nature.
The cardinality of primitive LFSRs of order $n$ over $\Fq$
is given by

\begin{equation} \label{NoLFSR}
\frac {\phi(q^{n}-1)}{n},
\end{equation}
where $\phi$ is Euler's totient function. Similarly
the number of irreducible LFSRs (whose characteristic polynomials are irreducible)
of order $n$ over a finite field $\Fq$ is given by
\begin{equation} \label{NoIrrLFSR}
 {\frac{1}{n} \displaystyle \sum_{d\mid n}
              \mu\left(d\right)q^{\frac{n}{d}}},
\end{equation}
where $\mu$ is the M\"{o}bius function.
\noindent
Zeng et.~al \cite{Zeng} considered a generalization of LFSR which they called as $\sigma$-LFSR.
These are word-oriented linear feedback shift registers, involving linear recurrence relation
over $\F_{2^m}$, with matrix coefficients coming from $M_m(\F_2)$. 
They also gave a conjectural formula for the number of primitive $\sigma$-LFSRs of order $n$ 
over $\F_{2^m}$ \cite{Zeng}. A further generalisation of $\sigma$-LFSR to $\Fqm$ was 
done by Ghorpade and Hasan \cite{GSM} who extended the conjecture formula over $\Fqm$. \cite{GSM} states
this number as
\begin{equation} \label{NoSigmaLFSR}
   \frac{\phi(q^{mn}-1)}{mn} q^{m(m-1)(n-1)}
        \displaystyle \prod_{i=1}^{m-1}(q^m-q^i).
\end{equation}
\noindent
Refer to \cite{GSM,GR,CT} for progress as well as complete proof of the conjecture.
Further refer to \cite{GR,CT,Ram} for the cardinality of irreducible $\sigma$-LFSRs which is given by 
\begin{equation}\label{NoIrrSigmaLFSR}
  {\displaystyle \frac{1}{mn}{q^{m(m-1)(n-1)}
   \displaystyle \prod_{i=1}^{m-1}(q^m-q^i)} \sum_{d\mid mn}
        \mu\left(d\right)q^{\frac{mn}{d}}}.
\end{equation}
\noindent
Here we focus on \emph{transformation shift registers} (TSRs) which are an extremely
important and useful subclass of $\sigma$-LFSRs. TSRs find their origin in a problem
posed by Bart Preneel \cite{BP} as a 
challenge to design fast and secure LFSRs which use the parallelism offered by the 
word operations of modern processors. The problem was addressed by 
the introduction of TSRs. Tsaban and Vishne \cite{TV} proved them to be faster 
and more efficient in software implementation than $\sigma$ LFSRs.
Refer to Dewar and Panario \cite{DP1,DP3} for further developments on the theory of TSRs.
Like $\sigma$ LFSRs, TSRs are also classified
as irreducible and primitive based on their characteristic polynomial. 
A study of irreducible TSRs was carried by Ram \cite{Ram} who considered
the problem of enumerating TSRs over a finite field and gave an explicit
formula for the number of irreducible TSRs of order two. The problem was
further investigated by Sartaj and Cohen \cite{SSDQ} who gave an asymptotic formula
for the number of irreducible TSRs in some special cases.So we see that some significant progress has 
been made on irreducible TSRs but the same cannot be said about primitive TSRs. In the context of stream ciphers, 
we are again basically interested in TSR sequences with maximum period i.e primitive TSRs.
Answers to questions regarding cardinality, existence, construction etc of primitive 
TSRs are challenging and still remain elusive. In our submission we concentrate on these aspects. 

We have made here an attempt to address the problem of existence and generation 
of primitive TSRs of order $n$ over $\F_{q^m}$ by establishing an equivalence
between primitive TSRs and primitive polynomials of special type. These primitive
polynomials serve as building blocks for primitive TSRs. We give a focussed search algorithm for 
generating primitive TSRs. A conjecture regarding the existence of these special type of primitive polynomials
has been proposed by us along with some experimental results in its support. We give an explicit proof for the existence of 
primitive TSRs of order $2$ over $\F_{2^m}$. Finally an attempt has been made by us to give the cardinality of primitive TSRs 
of order $2$ over $\F_{2^m}$ along with some bounds for primitive TSRs in general.

\section {Preliminaries} \label{sectsr}
\noindent
We will be using the following notations throughout the paper.
Let $\Fq$ denote the finite field with $q$ elements.
$\Fq[X]$ is the ring of polynomials with coefficients in $\Fq$. For 
every set $C$ let $|C|$ denote the cardinality of $C$. 
The set of all $d\times d$ matrices with entries in $\Fq$ is 
denoted by $M_d(\Fq)$. 

Throughout the paper we fix positive integers $m$ and $n$, 
and a vector space basis $\{\alpha_0, \dots, \alpha_{m-1}\}$
of ${\mathbb F}_{q^m}$ over $\mathbb F_q$. There exists a vector space isomorphism
from $\Fqm \longleftrightarrow \Fq^m$ such that $s \longmapsto \mathbf{s}$. 
here $\mathbf{s}$ denotes the corresponding co-ordinate vector $(s_0, \dots, s_{m-1})$ of $s$.
Elements of $\mathbb F_q^m$ may be thought of as row vectors
and so $\, \mathbf{s}C$ is a well-defined element of
$\mathbb F_q^m$ for any $\mathbf{s} \in \mathbb F_q^m$ and
$C\in M_m(\Fq)$.
We now recall from \cite{HPW} and \cite{Ram} some definitions and
results concerning transformation shift registers.

\begin{definition}\label{primitivepolynomial}
A polynomial $f(X) \in \Fq[X]$ of degree $n$ is said to be a primitive polynomial if its
root $\alpha$ generates the cyclic group $\F_{q^n}^*$, consisting of non zero
elements of $\F_{q^n}$.
\end{definition}

\begin{definition} \label{def:tsr} 
Let $c_0, c_1, \dots, c_{n-1} \in \Fq$ and $ A \in M_m(\mathbb F_q)$.
Given any $n$-tuple $(\mathbf{s}_0, \dots, \mathbf{s}_{n-1})$ of
elements of  $\mathbb F_{q^m}$, let $(\mathbf{s}_i)_{i=0}^{\infty}$
denote the infinite sequence of elements of ${\mathbb F}_{q^m}$
determined by the following linear recurrence relation:
\begin{eqnarray}\label{tsrdef}
{\mathbf{s}}_{i+n}={\mathbf{s}}_i(c_0A)+{\mathbf{s}}_{i+1}(c_1A)
          +\cdots +{\mathbf{s}}_{i+n-1}(c_{n-1}A)
          \quad i=0,1,\dots. \label{tsr}
\end{eqnarray}
The system (\ref {tsr}) is a \emph{transformation shift register}
(TSR) of order $n$ over $\mathbb F_{q^m}$, while the sequence
$(\mathbf{s}_i)_{i=0}^{\infty}$ is the \emph{sequence generated
by the TSR} (\ref{tsr}).
\begin{itemize}
\item The $n$-tuple $(\mathbf{s}_0,\mathbf{s}_1,\ldots, \mathbf{s}_{n-1})$ is the \emph{initial state} of the TSR.
\item The polynomial $I_mX^n -(c_{n-1}A)X^{n-1}- \cdots -(c_1A)X-(c_0A)$ with matrix coefficients is the \emph{tsr-polynomial}
		of the TSR (\ref{tsr}). Here $I_m$ denotes the $m\times m$ identity matrix over $\Fq$.
\item  The sequence $(\mathbf{s}_i)_{i=0}^{\infty}$ is \emph{ultimately periodic} if there are integers $r, n_0$ with
		$r\ge 1$ and $n_0\geq0$ such that $\mathbf{s}_{j+r}=\mathbf{s}_j$ for all $j \geq n_0$.
\item  The least positive integer $r$ with this property is the \emph{period} of $(\mathbf{s}_i)_{i=0}^{\infty}$
and the corresponding least nonnegative integer $n_0$ is the \emph{preperiod} of $(\mathbf{s}_i)_{i=0}^{\infty}$. The sequence
$(\mathbf{s}_i)_{i=0}^{\infty}$ is \emph{periodic} if its preperiod is $0$. 
\end{itemize}  
\end{definition}

We can associate a block companion matrix $T$ with a TSR definition given in (\ref{tsrdef}) as follows
\begin{equation} \label{typePP}
T =
\begin {pmatrix}
\mathbf{0} & \mathbf{0} & \mathbf{0} & . & . & \mathbf{0} & \mathbf{0} & c_0A\\
I_m & \mathbf{0} & \mathbf{0} & . & . & \mathbf{0} & \mathbf{0} & c_1A\\
. & . & . & . & . & . & . & .\\
. & . & . & . & . & . & . & .\\
\mathbf{0} & \mathbf{0} & \mathbf{0} & . & . & I_m & \mathbf{0} & c_{n-2}A\\
\mathbf{0} & \mathbf{0} & \mathbf{0} & . & . & \mathbf{0} & I_m & c_{n-1}A
\end {pmatrix},
\end{equation}
where $c_0, c_1, \dots , c_{n-1}\in \Fq$, $A\in M_m(\Fq)$ and 
$\mathbf{0}$
indicates the zero matrix in $M_m(\Fq)$. The set of all such
$(m,n)$-block companion matrices $T$ over $\Fq$ is denoted
by $\TSR(m,n,q)$. The block companion matrix \eqref{typePP} is the
state transition matrix for the TSR \eqref{tsr}. Indeed, the $k$-th
state $\mathbf{S}_k:=\left(\mathbf{s}_{k}, \mathbf{s}_{k+1}, \dots,
\mathbf{s}_{k+n-1}\right) \in \Fqm^n$ of the TSR (\ref{tsr}) is
obtained from the initial state
$\mathbf{S}_0:=\left(\mathbf{s}_{0}, \mathbf{s}_{1}, \dots,
\mathbf{s}_{n-1}\right) \in \Fqm^n$ by $\mathbf{S}_k = \mathbf{S}_0 T^k$,
for any $k\ge 0$.

Using a Laplace expansion or a suitable sequence of
elementary column operations, we conclude that if $T \in \TSR(m,n,q)$
is given by \eqref{typePP}, then $\det T = \pm \det (c_0A)$.
Consequently,
\begin{equation}
\label{nonsingP}
T \in \GLmn(\Fq) \Longleftrightarrow
    c_0 \neq 0 ~~\mbox{and}~~A\in \GLm(\Fq).
\end{equation}
where $\GLm(\Fq)$ is the general linear group of all $m \times m$
nonsingular matrices over $\Fq$. We denote here the intersection 
$\TSR(m,n,q) \cap \GLmn(\Fq)$ by $\TSR^*(m,n,q)$. Elements of $\TSR^*(m,n,q)$
are exactly  the state transition matrices of periodic TSRs
of order $n$ over $\Fqm$ \cite[Prop. 4]{HPW}. It follows from \eqref{typePP} that $T \in \TSR^*(m,n,q)$
iff $T$ is of the form 
\begin{equation} \label{typeP'}
\begin {pmatrix}
\mathbf{0} & \mathbf{0} & \mathbf{0} & . & . & \mathbf{0} & \mathbf{0} & B\\
I_m & \mathbf{0} & \mathbf{0} & . & . & \mathbf{0} & \mathbf{0} & c_1B\\
. & . & . & . & . & . & . & .\\
. & . & . & . & . & . & . & .\\
\mathbf{0} & \mathbf{0} & \mathbf{0} & . & . & I_m & \mathbf{0} & c_{n-2}B\\
\mathbf{0} & \mathbf{0} & \mathbf{0} & . & . & \mathbf{0} & I_m & c_{n-1}B
\end {pmatrix},
\end{equation}
\noindent
where $c_1, \dots, c_{n-1}\in \Fq$ and $B\in \GLm(\Fq)$. Henceforth
we deal with periodic TSRs only, that is, a TSR of the form (\ref{typeP'}).
\noindent
The map 
\begin{equation}\label{charmap}
\Psi : M_{mn}(\Fq) \longrightarrow \Fq[x] 
\end{equation}
defined by $\Psi(T) := det(XI_{mn} - T)$ will be referred to as the characteristic map.
The characteristic polynomial of $T$ is given by \cite[Lemma 1]{HPW}

\begin{equation}\label{CharT}
\Psi(T) = det(X^nI_m - g_T(X)B)
\end{equation}
where $g_T(X) = 1 + c_1X + c_2X + ...c_{n-1}X^{n-1} \in F_q[X]$. We see that $T$
is uniquely determined by $g_T(X)$ and $B$.
For every matrix A we denote by $\Psi_A(X)$ the characteristic polynomial of A.
It follows from ($\ref{CharT}$) that for $T \in \TSR^*(m,n,q)$

\begin{equation} \label{gencharT}
\Psi_T(X) = g_T(X)^m\Psi_B({\frac{X^n}{g_T(X)}})
\end{equation}
thus $f(X) \in \Psi(\TSR^*(m,n,q))$ iff $f(X)$ can be expressed in the form 

\begin{equation} \label{decom}
g(X)^mh(\frac{x^n}{g(X)})
\end{equation}
for some monic polynomial $h(X) \in \Fq[X]$ of degree $m$ with $h(0) \neq 0$ and 
$g(X) \in \Fq[X]$ of degree at most $n-1$ with $g(0)= 1$. When $f(X) \in \Psi(\TSR^*(m,n,q))$ is a primitive polynomial then the representation ($\ref{decom}$) is unique and is said to be $(m,n)$ decomposition of $f(X)$ \cite{Ram}.\\
\noindent

\section{Primitive TSRs} \label{TSRorder2}
\noindent
A TSR is \emph{primitive} 
if its characteristic polynomial is primitive.
The set of primitive TSRs is 
denoted by $\TSRP(m,n,q)$ and the set of primitive polynomials in 
$\Fq[X]$ of degree $d$ is denoted by $\mathcal P(d,q)$.
\\
\\
Then the \emph{characteristic map}
$$ \Psi: M_{mn}(\Fq) \longrightarrow \Fq[X]\quad \mbox{defined by}
        \quad\Psi(T):=\det(XI_{mn}-T),$$
if restricted to the set $\TSRP(m,n,q)$ 
yields the map
$$\Psi_{P}:\TSRP(m,n,q)\longrightarrow \mathcal P(mn,q).$$
It was noted in \cite{Ram} that the map $\Psi_{P}$ 
is not surjective in general. Denote the characteristic polynomial 
of $A \in M_{mn}(\Fq)$ by $\Psi_A(X)$
\\
\\
\\

\begin{lemma} \label{NoMatrices}
Let $ \eta: M_{m}(\Fq) \longrightarrow \Fq[X]$ be defined by
$\eta(A):=\det(XI_{m}-A)$. Then, for every $p(X) \in \mathcal
P(m,q)$, we have,
$$\left|\eta^{-1}\left(p(X)\right)\right| =
  \displaystyle \prod_{i=1}^{m-1}(q^m-q^i).$$
\end{lemma}
\noindent
\textbf{Proof:} \cite[Theorem 2]{Reiner}.

\subsection*{Primitive TSRs of odd order $n$ over $\Fqm, q \ge 3$ and $m \ge 2$}

\begin{theorem} \label{NoTSR}
The number of primitive TSRs of odd order $n$ over $\Fqm$ where $q \ge 3$ and $m \ge 2$
is given by 
$$\left|\TSRP(m,n,q)\right|
  = \left|\Psi_P\left(\TSRP(m,n,q)\right)\right|
    \displaystyle \prod_{i=1}^{m-1}(q^m-q^i).$$
\end{theorem}
\noindent
\textbf{Proof:}
Let us assume that $f(X) \in \Psi_P\left(\TSRP(m,n,q)\right)$ is the characteristic polynomial of $T \in \TSRP(m,n,q)$
i.e $\Psi_T(X) = f(X)$ then $f(X)$ can be uniquely expressed in the form (equation $12$)
\begin{equation} \label{irrform}
g(X)^m h\left(\frac{X^n}{g(X)}\right)
\end{equation}
$T \in \TSRP(m,n,q)$ $ \implies f(X)$ is primitive 
$\implies h(X)$ is primitive, (refer \cite{Ram}),
where $h(X)\in \Fq[X]$ of degree $m$
with $h(0)\neq 0$ and $g(X) \in \Fq[X]$ of
degree at most $n-1$ with $g(0)=1$.
Clearly $g_T(X)=g(X)$ and $\Psi_B(X)=h(X)$ refer (equation \ref{gencharT}). The number of such $T$
is equal to the number of possible values of $B$ with $\Psi_B(X)=h(X)$.
Since $h(X)$ is primitive, by Lemma \ref{NoMatrices}, the number
of such $B$ is $\displaystyle \prod_{i=1}^{m-1}(q^m-q^i)$ hence proved.
\\
\\
Let $P_{q}(m,n)$ denote the set of primitive polynomials of the form $X^{n} - \mu $g(X)$ $
$\in \Fqm[X]$ where $\mu$ is a primitive element of $\Fqm$, $g(X) \in \Fq[X]$ such that
$g(0) = 1$ and deg $g(X) \leq  (n-1)$.
\begin{theorem} \label{NoprimitiveTsr}
$|\TSRP(m,n,q)| = \frac{|P_q(m,n)|}{m}\frac{|\GLm(\Fq)|}{q^m-1}$ where $n$ is odd,
$q \ge 3$ and $m \ge 2$
\end{theorem}
\noindent
\textbf{Proof:} We will prove the above results along the lines of the proof of \cite[Theorem 6]{Ram}
\\
Define
$$\Omega_{q}(m,n) := \Psi_{P}(\TSRP(m,n,q)).$$
$\newline $  By theorem \ref{NoTSR}
$$|\TSRP(m,n,q)| = |\Omega_{q}(m,n)|\frac{\GLm(\Fq)}{q^m-1}.$$
\noindent
Define a map 

$$\Phi : P_{q}(m,n) \longrightarrow \Fqm[X] $$ by
$\newline$ 
$$\Phi(X^n - \mu g(X)) := \prod\limits_{i=0}^{m-1}(X^n - \mu^{q^i}g(X))).$$
\noindent
The product on the right is $(m,n)$ decomposable. Let $\beta$ be a root of $X^n - \mu g(X)$
in the extension field $\Fqmn$ then the minimal polynomial of $\beta$ over $\Fq$  is 
$\Phi(X^n - \mu g(X))$. Thus $\Phi(X^n - \mu g(X))$ is primitive in $\Fq[X]$. Since $\Omega_{q}(m,n)$
is precisely the set of primitive $(m,n)$ decomposable polynomials in $\Fq[X]$, it follows that 
$$\Phi(P_{q}(m,n)) \subseteq \Omega_{q}(m,n).$$ 
Claim is 
$$\Phi(P_{q}(m,n)) = \Omega_{q}(m,n).$$
\noindent
Let $f(x) \in \Omega_{q}(m,n) $. since $f$ is primitive, $f$ has a unique $(m,n)$ decomposition \cite[theorem 3]{Ram} say
$$ f(X) = {{g(X)}^m} h(\frac{X^n}{g(X)}).$$ 
$f(X)$ is primitive $\implies h(X)$ is primitive in $\Fq[X]$ and if $\mu$ is a root of $h(X)$ in $\Fqm$, then
$$ \Phi(X^n - \mu g(X)) = f(X).$$
Now $|\Phi^{-1}(f)| $ = $m$ for each $f$ $\in$ $\Omega_q(m,n)$ and therefore.
$$|\Omega_q(m,n)| = \frac{|P_q(m,n)|}{m}.$$

\subsection*{Primitive TSRs of order $n$ over $\F_{2^m}, m \ge 2$}
\begin{theorem}
The number of primitive TSRs of order $n$ over $\F_{2^m}$ where $m \ge 2$
is given by 
$$\left|\TSRP(m,n,2)\right|
  = \left|\Psi_P\left(\TSRP(m,n,2)\right)\right|
    \displaystyle \prod_{i=1}^{m-1}(2^m-2^i).$$
\end{theorem}
\noindent
\textbf{Proof:} Exactly along the line of the proof (Theorem \ref{NoTSR}) with $q=2$.

\begin{theorem}
$|\TSRP(m,n,2)| = \frac{|P_2(m,n)|}{m}\frac{|\GLm(\F_2)|}{2^m-1}$ where $n \ge 1$ 
and $m \ge 2$
\end{theorem}
\noindent
\textbf{Proof:} Exactly along the line of the proof (Theorem \ref{NoprimitiveTsr}) with $q=2$.

\section {Existence of primitive TSRs} \label{NoTSRoforder2}
\noindent
We denote by $n_{odd}$ whenever  $n$ is taken to be odd positive integer.\\

\noindent
Let $f(X) \in P_q(m,n_{odd})$ $\implies f(X) \in \F_{q^m}[X]$ and $f(X) =X^n-\mu g(X)$ where $\mu$ is a primitive element of $\Fqm$, $g(X) \in \Fq[X]$ such that
$g(0) = 1$ and deg $g(X) \leq  (n-1)$. Now Consider the reciprocal polynomial of $f(X)$ which is of the form $h(X) + \mu^{-1}$ where 
$h(X) \in \F_q[X], h(0)=0$ and $\mu^{-1} \in \F_{q^m}$ is a primitive element. Denote the reciprocal polynomials of $P_q(m,n_{odd})$
by $P(m,n_{odd},q)$. 

\medspace
\noindent
The existence of primitive TSRs of:
\begin{itemize}
\item odd order $n$ over $\Fqm$, $q \ge 3$ and $m \ge 2$ denoted by $\TSRP(m,n_{odd},q)$
\item any order $n \ge 2$ over $\F_{2^m}$, $m \ge 2$ denoted by $\TSRP(m,n,2)$
\end{itemize}
is directly connected to the problem of existence of primitive polynomials of the form $P_{q}(m,n_{odd}), q \ge 3$ and $m \ge 2$ whereas when $q=2$ it depends on primitive polynomials of the form $P_{2}(m,n)$ for any positive integer $m \ge 2,n \ge 2$.
\\
\\
Finally we have the following existence relation.\\
$$\TSRP(m,n_{odd},q) \iff P_q(m,n_{odd}) \iff P(m,n_{odd},q)$$ 
$$\TSRP(m,n,2) \iff P_2(m,n) \iff P(m,n,2)$$
Based on our experimental results (\ref{expresults}) we propose a conjecture regarding the existence of primitive polynomials $P(m,n,q)$ for any prime
$q$ and $m,n \ge 2$.

\begin{conjecture} \label{specprimpoly}
There exists a primitive polynomial $f(X)$ of degree $n$ over 
$\Fqm$ of the following form 
$$ f(X) = g(X) + \lambda, $$ 
$\forall \ m, n \ge 2$ and $\forall q$, where $g(X) \in \Fq[X]$ such that
$g(0) = 0$ and $\lambda$ is a primitive element in $\Fqm$.\\
\end{conjecture}

\noindent
Denote the Galois group of automorphisms of $\F_{q^m}$ over $\F_q$ by $Gal(\F_{q^m}/\F_q)$ then a useful and alternate form of the conjecture \ref{specprimpoly} is as follows:
\begin{conjecture} \label{specprimpolyalt}
For all $m,n$ there exist polynomials $f(X), g(X) \in F_q[X]$ of degrees $m$ and $n$ respectively with $f(X)$  primitive and $g(0)=0$ such that $f(g(x)) \in \Fq[X]$ is primitive of degree $mn$.
\end{conjecture}
\noindent
\textbf{Proof:} Suppose $f(X)=g(X)+\lambda,$ as described in conjecture (\ref{specprimpoly}).\\

$\implies \prod\limits_{\sigma \in Gal(\F_{q^m}/\F_q)}^{} \sigma(f(X)) \in \Fq[X]$ is primitive of degree $mn$\\
Now $h(X)=\prod\limits_{i=0}^{m-1}(X+\lambda^{q^i}) \in \Fq[X]$ is primitive of degree $m$\\
but $h(g(X))= \prod\limits_{\sigma \in Gal(\F_{q^m}/\F_q)}^{} \sigma(f(X))$ is primitive of degree $mn$.\\
\noindent
Therefore we have $h(X) \in \Fq[X]$ primitive polynomial of degree $m$ and $g(X) \in \Fq[X]$ of degree $n$ such that $g(0)=0.$\\
Conversely\\ let $f(X), g(X) \in \Fq[X]$ be as given in conjecture (\ref{specprimpolyalt}) such that $f(g(X)) \in \Fq[X]$ is a primitive polynomial of degree $mn$.\\
Now $ f(X) = \prod\limits_{i=0}^{m-1}(X+\lambda^{q^i})$, $\lambda^{q^i}$ are primitive roots of $f(X)$ in $\Fqm$ for $ i \in \{0, \dots m-1\}$. Therefore
$f(g(X))=\prod\limits_{i=0}^{m-1}(g(X)+\lambda^{q^i})$ is primitive $\implies g(X)+\lambda^{q^i}$ is primitive $\forall i \in \{0, \dots m-1\}$. Hence $h(X)=g(X)+\lambda$ is a primitive polynomial in $\Fqm$.\\

\noindent
We now give a search algorithm for generating primitive TSRs of odd order $n$ over $\Fqm$.
\section{Search algorithm for primitive TSRs of odd order $n$ over $\F_{q^m}, q \ge 3$} \label{psearchforoddn}

\begin{itemize}

\item[{\rm step 1.}] Pick a primitive polynomial $f(X)$ of degree $m$ over $\F_q$.
\item[{\rm step 2.}] Pick a polynomial $g(X)$ of odd degree $n$ in $F_q$ such that $g(0) = 0$.
\item[{\rm step 3.}] Check if $f(g(X))$ primitive over $F_q$.
\item[{\rm step 4.}] If primitive, proceed to step $5$ else repeat step $1$.
\item[{\rm step 5.}] Take $k(X)=g(X) + \alpha$ such that $f(\alpha) = 0$. Compute the reciprocal polynomial of $k(X)$ given by $X^n + \lambda(X^ng(\frac{1}{X}))$ where $\lambda = \alpha^{-1}$. Therefore\\ reciprocal($k(X)$)
					 = $X^n + \lambda(c_{n-1}X^{n-1} + c_{n-2}X^{n-2}....c_1X + 1)=X^n+{\lambda}L(X)$
\item[{\rm Step 6.}] Compute the minimal polynomial, say $h(X)$, of $\lambda$ in $\F_q[X]$. It is primitive.
\item[{\rm Step 7.}] Compute matrix $A$ in $GL_m(\F_q)$ whose characteristic polynomial is $h(X)$. 
\item[{\rm Step 8.}] The characteristic polynomial of TSR $T$ is 
					 $$\prod\limits_{\sigma \in Gal(\F_{q^m}/\F_q)}^{} \sigma(X^n + \lambda L(X)).$$
					 It is primitive in $\F_q[X]$.
\item[{\rm Step 9.}] $T$ is given by 
					\begin{equation} \label{ConstTSR}
					T =
					\begin {pmatrix}
					\mathbf{0} & \mathbf{0} & \mathbf{0} & . & . & \mathbf{0} & \mathbf{0} & A\\
					I_m & \mathbf{0} & \mathbf{0} & . & . & \mathbf{0} & \mathbf{0} & c_1A\\
					. & . & . & . & . & . & . & .\\
					. & . & . & . & . & . & . & .\\
					\mathbf{0} & \mathbf{0} & \mathbf{0} & . & . & I_m & \mathbf{0} & c_{n-2}A\\
					\mathbf{0} & \mathbf{0} & \mathbf{0} & . & . & \mathbf{0} & I_m & c_{n-1}A
					\end {pmatrix}.
					\end{equation}
					
\end{itemize}
\noindent
\noindent
The above algorithm can be exactly used for generating primitive TSRs of any order $n$ over $\F_{2^m}$.
However for clarity we restate the algorithm.

\section{Search algorithm for primitive TSRs of order $n$ over $\F_{2^m}$} \label{psearchforalln}
\begin{itemize}
\item[{\rm step 1.}] Pick a primitive polynomial $f(X)$ of degree $m$ over $\F_2$.
\item[{\rm step 2.}] Pick a polynomial $g(X)$ of degree $n$ in $F_2$ such that $g(0) = 0$.
\item[{\rm step 3.}] Check if $f(g(X))$ primitive over $F_2$.
\item[{\rm step 4.}] If primitive, proceed to step $5$ else repeat step $1$.
\item[{\rm step 5.}] Take $k(X)=g(X) + \alpha$ such that $f(\alpha) = 0$. Compute the reciprocal polynomial of $k(X)$ given by $X^n + \lambda(X^ng(\frac{1}{X}))$ where $\lambda = \alpha^{-1}$. Therefore\\ reciprocal($k(X)$)
					 = $X^n + \lambda(c_{n-1}X^{n-1} + c_{n-2}X^{n-2}....c_1X + 1)=X^n+{\lambda}L(X)$
\item[{\rm Step 6.}] Compute the minimal polynomial, say $h(X)$, of $\lambda$ in $\F_2[X]$. It is primitive.
\item[{\rm Step 7.}] Compute matrix $A$ in $GL_m(\F_2)$ whose characteristic polynomial is $h(X)$. 
\item[{\rm Step 8.}] The characteristic polynomial of TSR $T$ is 
					 $$\prod\limits_{\sigma \in Gal(\F_{2^m}/\F_2)}^{} \sigma(X^n + \lambda L(X)).$$
					 It is primitive in $\F_2[X]$.
\item[{\rm Step 9.}] $T$ is given by 
					\begin{equation} \label{ConstTSR}
					T =
					\begin {pmatrix}
					\mathbf{0} & \mathbf{0} & \mathbf{0} & . & . & \mathbf{0} & \mathbf{0} & A\\
					I_m & \mathbf{0} & \mathbf{0} & . & . & \mathbf{0} & \mathbf{0} & c_1A\\
					. & . & . & . & . & . & . & .\\
					. & . & . & . & . & . & . & .\\
					\mathbf{0} & \mathbf{0} & \mathbf{0} & . & . & I_m & \mathbf{0} & c_{n-2}A\\
					\mathbf{0} & \mathbf{0} & \mathbf{0} & . & . & \mathbf{0} & I_m & c_{n-1}A
					\end {pmatrix}.
					\end{equation}
\end{itemize}

\noindent
We now give an explicit proof for the existence of primitive TSRs of order $2$ over $\Fqm$.

\begin{theorem} \label{Traceoneprimitiveelement}
\cite{OM} There exists a primitive quadratic polynomial of trace $1$ over $\F_{2^{m}}$.
\end{theorem} 

\begin{corollary} \label{specprim}
There exists a primitive quadratic polynomial of the form $X^2+ \lambda X+\lambda$ in $\F_{2^m}[X]$, $\F^*_{2^m}=<\lambda>$ $\forall m \ge 1$.
\end{corollary}

\begin{theorem} \label{primtiveTsrExistence}
There exists a primitive TSRs of order $2$ over $\F_{2^m}$ for all m.
\end{theorem}

Proof:- Using corollary (\ref{specprim}) consider a primitive polynomial of the form
$X^2+ \lambda X+\lambda$. Now consider a map 
$$\F_{2^{m}}[X] \longleftarrow \F_{2}[X]$$
$$f(X) \longmapsto \prod\limits_{\sigma \in Gal(\F_{2^m}/\F_2)}^{} \sigma(f(X))$$
$$X^2+ \lambda X+\lambda \longmapsto \prod\limits_{\sigma \in Gal(\F_{2^m}/\F_2)}^{} \sigma(X^2+ \lambda (X+1))$$
but
$$\sigma(X^2+ \lambda (X+1)) = (X+1)\{\frac{X^2}{(X+1)}+\sigma(\lambda)\}$$
$$\prod\limits_{\sigma \in Gal(\F_{2^m}/\F_2)}^{} \sigma(X^2+ \lambda (X+1)) =  g(X)h(\frac{X^2}{g(X)})$$
\\
where
$$g(X) = (X+1) \ and \ h(X) = \prod\limits_{\sigma \in Gal(\F_{2^m}/\F_2)}^{}(X-\sigma(\lambda))$$
\\
Now $g(X)h(\frac{X^2}{g(X)})$ is a primitive polynomial of degree $2m$ which gives a primitive 
TSR of order $2$ over field $F_{2^m}$.

\section{Cardinality of $P_2(m,2)$}

\noindent
We now consider the cardinality of primitive TSRs of order $n$ over $\Fqm$ for trivial values of $m$ and $n$.\\
\noindent
The case $n=1$ follows immediately from \cite[Theorem 7.1]{GSM}.
In this case, the number of primitive TSRs of order one over
$\Fqm$ is given by
$$
  \frac{\left|\GLm(\Fq)\right|}{(q^m-1)}\frac {\phi(q^m-1)}{m}.
$$
\noindent
The case $m=1$ is trivial and in this case, the number of primitive
TSRs of order $n$ is given by 
$$\frac {\phi(q^{n}-1)}{n}.$$
\noindent
However, for general values of $m$ and $n$, the enumeration of
primitive TSRs does not seem to be an easy problem and remains open. We attempt to derive the cardinality of
primitive TSRs of order $2$ over $\F_{2^m}$. Consider
$$P_2(m,2) := \{f(X) = X^2 - \alpha(X + 1) : f(X) \ \mbox{primitive} , f(X) \in \F_{2^m}[X], \F_{2^m}^{*} = \ <\alpha>\}.$$
$$|P_2(m,2)| = |\{X^2 + X + \alpha : \F_{2^m}^{*} = \ <\alpha>\} |.$$
\noindent
Consider the field $\F_{2^{2m}}$ generated by a primitive polynomial $f(X) \in \F_2[X]$ of degree $2m$ and
having a primitive root $\alpha \in \F_{2^{2m}}$. Consider the set 
$$A = \{i : (i,2^{2m}-1) = 1, 1 \leq i \leq 2^{2m}-1 \}.$$
Then the set
$$B = \{\alpha^i: i \in A \}, $$ gives the primitive elements in $\F_{2^{2m}}^*$. 

\begin{definition} \label{cosetclass}
\cite{GG} A cyclotomic coset $C_{i_{s}}$ of an element $i_s \in D$ modulo $2^n - 1$ with respect to $2$ is defined to be 
$$C_{i_{s}} = \{i_s, i_s.2, i_s.2^2, i_s.2^3....i_s.2^{n_{s}-2}, i_s.2^{n_{s}-1}\}$$
where $n_s$ is the smallest positive integer such that $i_s \equiv i_s.2^{n_s} (mod \ 2^n-1)$. The subscript 
$i_s$ is chosen as the smallest integer in $C_s$, and $i_s$ is called the coset leader of $C_s$. 
\end{definition}
\noindent
Decompose A in cyclotomic cosets with respect to 2 modulo $2^{2m}-1$. Let $C_{i_{1}}$ 
denote the cyclotomic coset containing $i_1 \in A$. Corresponding to $C_{i_1}$, $\alpha_{C_{i_1}}$ is the set 
containing $\alpha^{i_1} \in B$ and all its conjugates over $\F_2$. We call this as conjugate class of $C_{i_1}$
$$C_{i_{1}} = \{i_1, i_1.2, i_1.2^2, i_1.2^3....i_1.2^{2m-2}, i_1.2^{2m-1}\},$$
$$\alpha_{C_{i_{1}}} = \{\alpha^{i_1}, \alpha^{{i_1}.2}, \alpha^{{i_1}.2^3}...\alpha^{{i_1}.2^{2m-2}},\alpha^{{i_1}.2^{2m-1}}\}.$$
where $(i_1,2^{2m}-1) = 1$, $\alpha^{2^{2m}-1} = 1$ and $|C_{i_1}| = 2m$. 
Let $C_{i_1}, C_{i_2},...C_{i_n}$ be the cyclotomic coset classes of $A$ then
\begin{eqnarray*}
A &=& \cup^{j=n}_{j=1}C_{i_j},\\
B &=& \cup^{j=n}_{j=1}\alpha_{C_{i_{j}}}.
\end{eqnarray*}
\noindent
For each conjugate class consisting of elements of B form polynomials
\begin{eqnarray*}
X^2&-&(\alpha^j + \alpha^{j.2^m})X + \alpha^j.\alpha^{j.2^m},\\
X^2&-&(\alpha^{j.2} + \alpha^{j.{2^{m+1}}})X + \alpha^{j.2}.\alpha^{j.{2^{m+1}}},\\
X^2&-&(\alpha^{j.2^2} + \alpha^{j.{2^{m+2}}})X + \alpha^{j.2^2}.\alpha^{j.{2^{m+2}}},\\
.\\
.\\
.\\
X^2&-&(\alpha^{j.2^{m-1}} + \alpha^{j.{2^{2m-1}}})X + \alpha^{j.2^{m-1}}.\alpha^{j.{2^{2m-1}}}.
\end{eqnarray*}
\noindent
This gives the $m$ polynomials given by the conjugate class $\alpha_{C_j}$ corresponding to cyclotomic class $C_{j}$. Here we have trace and norm of polynomials as
$$ N = \{\alpha^{j+j.2^m}, \alpha^{({j + j.2^m}).2}, \alpha^{({j + j.2^m}).2^2}...\alpha^{({j + j.2^m}).2^{(m-1)}}\}, $$
$$ T = \{{(\alpha^j + \alpha^{j.2^m})}, {(\alpha^j + \alpha^{j.2^m})^2}, {(\alpha^j + \alpha^{j.2^m})^{2^2}}...{(\alpha^j + \alpha^{j.2^m})^{2^{(m-1)}}}  \}.$$
\\
\noindent 
In one of the conjugate class of elements of $B$ we will get primitive elements with trace one (\ref{Traceoneprimitiveelement}).
Consequently trace of all the quadratic polynomials formed from that class will be $1$. Counting such conjugate classes will
give us the number of primitive quadratic polynomials with trace $1$. 
Since each conjugate class gives, $m$ degree $2$, polynomials. Therefore 
total number of primitive quadratic polynomials over $\F_{2^m}$ with trace $1$ will be multiple of $m$ 

$$|P_2(m,2)| = rm,$$
where $r$ is the number of conjugate classes with trace $1$.
Hence the number of primitive TSRs of order $2$ over $\F_{2^{m}}$ is
$$|\TSRP(m,2,2)| = \frac{|P_2(m,2)|}{m}\frac{|\GLm(\F_2)|}{2^m-1},$$
$$|\TSRP(m,2,2)| = \frac{r*|\GLm(\F_2)|}{2^m-1}, $$

$r \le \frac{\phi(2^m-1)}{m}$. Here we give some values of $r$ obtained experimentally using sage.

\begin{table}[htbp]
\begin{center}
\small
\label{mid1}
{
\begin{tabular} {|p{.5cm}|p{.5cm}| p{4.5cm}|}
\hline
m & r & No of primitive quadratic polynomials over $\F_{2^m}$ with trace one. $|P_2(m,2)|$ \\
\hline
2 & 1 & 2 \\
3 & 1 & 3 \\
4 & 1 & 4 \\
5 & 2 & 10 \\
6 & 3 & 18 \\
7 & 6 & 42 \\
8 & 7 & 56\\
9 & 16 & 144\\
10 & 25 & 250\\
11 & 57 & 627\\
12 & 68 & 816\\
\hline
\end{tabular}
}
\end{center}
\end{table}

\section{Bounds on the number of primitive TSRs in special cases}

%

\begin{theorem} \label{numbound}

$$|\TSRP(m,n,q)| \le (q^{n-1} - 1)\frac{\phi(q^m - 1)}{m}\frac{\GLm(\Fq)}{q^m - 1},$$ where $n$ is odd.
In particular, TSRs of order $2$ over $\F_{2^m}$ are bounded by,
$$|\TSRP(m,2,2)| \le (2^{n-1} - 1)\frac{\phi(2^m - 1)}{m}\frac{\GLm(\F_2)}{2^m - 1}.$$
\end{theorem}

Proof: Consider $n$ to be odd then from (\ref{NoprimitiveTsr}) we have, 
$$|\TSRP(m,n,q)| = \frac{|P_q(m,n)|}{m}\frac{|\GLm(\Fq)|}{q^m-1}.$$
Consider all primitive polynomials of deg $n$ in $\Fqm[x]$ of the form
$$ \U = {\{g(X) + \lambda : g(X) \in \F_q[X], \ g(0) = 0, \ \lambda \ \mbox{primitive in} \ \Fqm \}},$$
$$ |P_{q}(m,n)| = |\U|,$$
but $|\U| \le (q^{n-1} - 1)\phi(q^m - 1). $
\newline
Similarly for the case when $q = 2$.
\noindent

\section{Some experimental verification for the proposed conjecture}\label{expresults}
\noindent
The conjecture (\ref{specprimpoly}) is always true for $n = 2$ and $q = 2$ since we always have primitive elements in $\F_{2^{2m}}$ with 
trace $1$ over $\F_{2^m}$ for any $m \ge 2$. Therefore in $\F_{2^m}[X]$ we always have primitive polynomials $f(X)$ of degree $2$
of the form 
$$f(X) = g(X) + \lambda,$$
where $\lambda$ is a primitive element in $\F_{2^m}$ and $g(X) = X^2 + X$. However based on our experimental results we feel
that the conjecture is always true. We give below some results in support of our claim. We recall from the conjecture (\ref{specprimpoly}) that $P(m,n,q)$ denotes the primitive polynomials of degree $n$ over $\Fqm$ of the form $g(X) + \lambda$ where $g(X) \in \Fq[X]$ is of degree $n$ such that $g(0)=0$ and $\lambda$ is a primitive element of $\Fqm$.
\\
\begin{table}[H]
\begin{center}
\tiny
\small
\caption{Primitive polynomials of the type $P(2,3,q)$ i.e of degree $3$ over $\F_{q^2}$ where $\F_{q^2}^* = <a>$ } 
\label{mid}
{
\begin{tabular} {|p{1.0cm}|p{6.0cm}|}
\hline
$q$ & primitive polynomials\\
\hline
& \\
2 & $x^3 + x^2 + x + a$\\
  & $x^3 + x^2 + x + a+1$\\
  & \\
3 & $x^3 + x^2 + x + a$\\
  & $x^3 + x^2 + x + 2a+1$\\
  & \\
5 & $x^3 + x^2 + x + 3a$\\
  & $x^3 + x^2 + x + 2a+3$\\
  & \\
7 & $x^3 + x^2 + x + 3a+1$\\
  & $x^3 + x^2 + x + 3a+3$\\
  & $x^3 + x^2 + x + 4a+4$\\
  & $x^3 + x^2 + x + 4a+6$\\
  & \\
11& $x^3 + x^2 + 7a+1$\\
  & $x^3 + x^2 + a+10$\\
  & $x^3 + x^2 + a+7$\\
  & $x^3 + x^2 + 4a+7$\\
  & $x^3 + x^2 + 3a$\\
  & $x^3 + x^2 + 8a+1$\\
  & \\
\hline

\end{tabular}
}
\end{center}
\end{table}

\begin{table}[H]
\begin{center}
\tiny
\small
\caption{Primitive polynomials of the type $P(m,3,2)$ i.e of degree $3$ over $\F_{2^m}$ where $\F_{2^m}^* = <a>$ } 
\label{mid1}
{
\begin{tabular} {|p{1.0cm}|p{6.0cm}|}
\hline
$m$ & primitive polynomials\\
\hline
& \\
3 & $x^3 + x^2 + a$\\
  & $x^3 + x^2 + a^2$\\
  & $x^3 + x^2 + a^2+a$\\
  & \\
4 & $x^3 + x^2 + a^3+a+1$\\
  & $x^3 + x^2 + a^3+a^2+a$\\
  & $x^3 + x^2 + a^3+a^2+1$\\
  & $x^3 + x^2 + a^3+1$\\
  & \\
5 & $x^3 + x^2 + a^4+a^2$\\
  & $x^3 + x^2 + a^2+a+1$\\
  & $x^3 + x^2 + a^4+a^3+a^2$\\
  & $x^3 + x^2 + a^4+a^3+a^2+1$\\
  & $x^3 + x^2 + a^2+a$\\
  & $x^3 + x^2 + a^4+a^3$\\
  & $x^3 + x^2 + a^4+a^2+1$\\
  & $x^3 + x^2 + a^4+a^3+1$\\
  & $x^3 + x^2 + a^4+a^2+a+1$\\
  & $x^3 + x^2 + a^4+a^2+a$\\
  & \\
6 & $x^3 + x^2 + a^4+a^3+1$\\
  & $x^3 + x^2 + a^5+a^4+a^3+a$\\
  & $x^3 + x^2 + a^5+a^3+a^2$\\
  & $x^3 + x^2 + a^4+a^3$\\
  & $x^3 + x^2 + a^5+a^4+a^3+a+1$\\
  & $x^3 + x^2 + a^5+a^3+a^2+1$\\
  & \\
\hline

\end{tabular}
}
\end{center}
\end{table}

\begin{table}[H]
\begin{center}
\tiny
\small
\caption{Primitive polynomials of the type $P(2,n,2)$ i.e. of degree $n$ over $\F_{2^2}$ where $\F_{2^2}^* = <a>$ i.e } 
\label{mid2}
{
\begin{tabular} {|p{1.0cm}|p{6.0cm}|}
\hline
$n$ & primitive polynomials\\
\hline
& \\
4 & $x^4 + x^3 + x^2 + a$\\
  & $x^4 + x^3 + x^2 + a+1$\\
  & \\
5 & $x^4 + x^3 + x^2 + x + a$\\
  & $x^4 + x^3 + x^2 + x + a+1$\\
  & \\
6 & $x^6 + x^5 + x + a$\\
  & $x^6 + x^5 + x + a+1$\\
  & \\
7 & $x^7 + x^6 + x^5 + a$\\
  & $x^7 + x^6 + x^5 + a+1$\\
  & $x^7 + x^6 + x^4 + a$\\
  & $x^7 + x^6 + x^4 + a+1$\\
  & $x^7 + x^4 + x^3 + a$\\
  & $x^7 + x^4 + x^3 + a+1$\\
  & $x^7 + x^6 + x^4 + x^3 + a$\\
  & $x^7 + x^6 + x^4 + x^3 + a+1$\\
  & $x^7 + x^6 + x^2 + a$\\
  & $x^7 + x^6 + x^2 + a+1$\\
  & $x^7 + x^5 + x^4 +x^2 + a$\\
  & $x^7 + x^5 + x^4 + x^2 + a+1$\\
  & $x^7 + x^6 + x^5 + x^4 + x^3 + x^2 + a$\\
  & $x^7 + x^6 + x^5 + x^4 + x^3 + x^2  + a+1$\\
  & $x^7 + x^5 + x^4 + x + a$\\
  & $x^7 + x^5 + x^4 + x + a+1$\\
  & $x^7 + x^6 + x^5 + x^4 + x + a$\\
  & $x^7 + x^6 + x^5 + x^4 + x + a+1$\\
  & $x^7 + x^3 + x + a$\\
  & $x^7 + x^3 + x + a+1$\\
  & $x^7 + x^5 + x^4 + x^3 + x + a$\\
  & $x^7 + x^5 + x^4 + x^3 + x + a+1$\\
  & $x^7 + x^5 + x^3 + x^2 + x + a$\\
  & $x^7 + x^5 + x^3 + x^2 + x + a+1$\\
  & $x^7 + x^6 + x^5 + x^3 + x^2 + x + a$\\
  & $x^7 + x^6 + x^5 + x^3 + x^2 + x + a+1$\\
  & $x^7 + x^6 + x^4 + x^3 + x^2 + x + a$\\
  & $x^7 + x^6 + x^4 + x^3 + x^2 + x + a+1$\\
  & \\
\hline

\end{tabular}
}
\end{center}
\end{table}

\begin{table}[H]
\begin{center}
\tiny
\small
\caption{Primitive polynomials of the type $P(m,3,3)$ i.e.of degree $3$ over $\F_{3^m}$ where $\F_{3^m}^* = <a>$. i.e } 
\label{mid3}
{
\begin{tabular} {|p{1.0cm}|p{6.0 cm}|}
\hline
$m$ & primitive polynomials\\
\hline
& \\
3 & $x^3 + x^2 + a$\\
  & $x^3 + x^2 + a+2$\\
  & $x^3 + x^2 + a^2+2a+2$\\
  & $x^3 + x^2 + a+1$\\
  & $x^3 + x^2 + a^2+a+2$\\
  & $x^3 + x^2 + 2a^2+a$\\
  & $x^3 + x^2 + a^2+1$\\
  & $x^3 + x^2 + 2a^2+2a$\\
  & $x^3 + x^2 + 2a^2+1$\\
  & \\
4 & $x^3 + x + a$\\
  & $x^3 + x + a^3$\\
  & $x^3 + x + 2a^3+a^2+a+1$\\
  & $x^3 + x + a^3+a^2+2a$\\
  & $x^3 + x + a^3+a+2$\\
  & $x^3 + x + 2a^3+a^2+2a$\\
  & $x^3 + x + 2a^3+2a$\\
  & $x^3 + x + a^3+2a+2$\\
  & $x^3 + x + 2a^2+a+1$\\
  & $x^3 + x + a^3+2a^2+1$\\
  & $x^3 + x + a^2+a$\\
  & $x^3 + x + 2a^3+a^2+2a+2$\\
  & $x^3 + x + 2a$\\
  & $x^3 + x + 2a^3$\\
  & $x^3 + x + a^3+2a^2+2a+2$\\
  & $x^3 + x + a^3+2a^2+a$\\
  & $x^3 + x + 2a^3+2a+1$\\
  & $x^3 + x + a^3+2a^2+a$\\
  & $x^3 + x + a^3+a$\\
  & $x^3 + x + 2a^3+a+1$\\
  & $x^3 + x + a^2+2a+2$\\
  & $x^3 + x + 2a^3+a^2+2$\\
  & $x^3 + x + 2a^2+2$\\
  & $x^3 + x + a^3+2a^2+a+1$\\
  & $x^3 + x + 2a^3+a^2+a+1$\\
  & $x^3 + x + 2a^3+2a^2+a+1$\\
  & $x^3 + x + a^3+2a+2$\\
  & $x^3 + x + 2a^2+a+1$\\
  & $x^3 + x + a^2+a$\\
  & $x^3 + x + 2a^3+1$\\
  & $x^3 + x + 2a+1$\\
  & $x^3 + x + 2a^3+a^2$\\
  & $x^3 + x + a^3+2a^2+2a+2$\\
  & $x^3 + x + 2a^3+a+1$\\
  & $x^3 + x + a^2+2a+2$\\
  & $x^3 + x + 2a^2+2a$\\
  & \\
\hline

\end{tabular}
}
\end{center}
\end{table}
\noindent
Based on some specific observations we propose the existence of trace $1$ trinomials of the special form over $\F_{q^2}$ for all prime $q$.

\begin{conjecture} \label{cubicpoly}
There always exist primitive polynomials of the type $P(2,3,q)$ which have the form $x^3+x^2+x+\alpha$ over $\F_{q^2} \  \forall \ q$ where $\F_{q^2}^* = <a>$.
\end{conjecture}

\begin{table}[H]
\begin{center}
\tiny
\small
\caption{Primitive polynomials of degree $3$ over $\F_{q^2}$ where $\F_{q^2}^* = <a>$} 
\label{mid1}
{
\begin{tabular} {|p{1.5cm}|p{4.5cm}|}
\hline
$q$ & primitive polynomials\\
\hline
& \\
2 & $x^3 + x^2 + x + a$\\
  & $x^3 + x^2 + x + a+1$\\
  & \\
3 & $x^3 + x^2 + x + a$\\
  & $x^3 + x^2 + x + 2a+1$\\
  & \\
5 & $x^3 + x^2 + x + 3a$\\
  & $x^3 + x^2 + x + 2a+3$\\
  & \\
7 & $x^3 + x^2 + x + 3a+1$\\
  & $x^3 + x^2 + x + 3a+3$\\
  & $x^3 + x^2 + x + 4a+4$\\
  & $x^3 + x^2 + x + 4a+6$\\
  & \\
11& $x^3 + x^2 + x + 9a+2$\\
  & $x^3 + x^2 + x + 9a+6$\\
  & $x^3 + x^2 + x + 6a+5$\\
  & $x^3 + x^2 + x + 5a$\\
  & $x^3 + x^2 + x + 6a+4$\\
  & $x^3 + x^2 + x + 6a+9$\\
  & $x^3 + x^2 + x + 2a+9$\\
  & $x^3 + x^2 + x + 2a+5$\\
  & $x^3 + x^2 + x + 5a+6$\\
  & $x^3 + x^2 + x + 6a$\\
  & $x^3 + x^2 + x + 5a+7$\\
  & $x^3 + x^2 + x + 5a+2$\\
  & \\
13& $x^3 + x^2 + x + a$\\
  & $x^3 + x^2 + x + 12a+6$\\
  & $x^3 + x^2 + x + 10a+9$\\
  & $x^3 + x^2 + x + 12a+1$\\
  & $x^3 + x^2 + x + 11a+9$\\
  & $x^3 + x^2 + x + 7a+5$\\
  & $x^3 + x^2 + x + 9a+5$\\
  & $x^3 + x^2 + x + 10a+11$\\
  & $x^3 + x^2 + x + 2a+7$\\
  & $x^3 + x^2 + x + a+5$\\
  & $x^3 + x^2 + x + 4a+1$\\
  & $x^3 + x^2 + x + 9a$\\
  & \\
\hline

\end{tabular}
}
\end{center}
\end{table}

\section{Discussions} 
\noindent
Primitive TSRs are very important for generating efficient word oriented stream ciphers. We have here dealt 
with the question of existence of primitive TSRs of order $2$ over $F_{2^m}$ and attempted to derive a formula for 
their cardinality. However, we saw that the problem of computing the number of primitive TSRs of order $2$ 
is related to computing the number of primitive elements in $\F_{2^{2m}}$ with trace $1$ over $F_{2^m}$. 
Also a general construction algorithm for finding primitive TSRs of order $2$ over $F_{2^m}$ is related to the construction 
algorithm for finding primitive elements in $\F_{2^{2m}}$ with trace $1$ over $F_{2^m}$. As far as the question of existence of primitive TSRs
of odd order $n$ over $\Fqm$, where $q \ge 2$, is concerned, we have proposed a conjecture (\ref{specprimpoly}) regarding the 
existence of primitive polynomials of special type $P(m,n,q)$. We propose the following questions for further
study.

\begin{itemize}
\item[{\rm 1.}] Computing the number of conjugate classes of primitive elements in $\F_{2^{2m}}$ with trace $1$ over $F_{2^m}$.
\item[{\rm 2.}] Construction algorithm for finding primitive elements in $\F_{2^{2m}}$ with trace $1$ over $F_{2^m}$.
\item[{\rm 3.}] Existence of primitive polynomials of the form $P(m,n,q)$ and $P(m,n,2)$.
\end{itemize}

\bibliography{etsr}{}
\bibliographystyle{plain}

\end{document}